\documentclass[leqno,12pt]{article}

\usepackage{amsmath,amsthm,amsfonts,amssymb}

\newtheorem{thm}{Theorem}
\newtheorem{lem}{Lemma}

\newcommand{\m}{\mathcal}
\newcommand{\e}{\varepsilon}

\title{\bf On the distribution of $r$-tuples of squarefree numbers in short intervals }

\author{ \bf D. I. Tolev  }

\date{ }

\begin{document}

\maketitle

\begin{abstract}
We cosider the number of $r$-tuples of squarefree numbers in a 
short interval. 
We prove that it cannot be 
much bigger than the expected value and we also 
estabish an asymptotic formula
if the interval is not very short.

\end{abstract}

\section{Introduction and statement of the results.}

\indent

For $x \ge 1$ and $ h \ge 1 $ we define
\[
  Q(x)= \sum_{ n \le x} \mu^2(n) , \qquad
  Q(x, h) = Q(x+h) - Q(x), 
\]
where $\mu(n) $ is the M\"obius function. It is conjectured that 
for any  $ \e >0 $ there exists $ x_0(\e) $ such that 
$ Q(x, x^{\e}) > 0 $ whenever $ x \ge x_0(\e) $.
A conditional proof of this hypothesis (under the {\it ABC}-conjecture)
was found in 1998 by Granville \cite{G}. An 
unconditional proof is not known at present, 
but many approximations were established during the last decades.
The strongest of them is due to Filasaeta and Trifonov \cite{FT}
(an information about the earlier work on this problem 
is also available there).
In 1992, using clever elementary arguments, they proved the following:
\begin{thm}[Filaseta, Trifonov] \label{thmFT}
There exists a constant $c>0$ such that if $x$ is sufficiently large
and $ h = c \, x^{1/5} \log x $,
then $ Q(x, h) > 0 $.
\end{thm}
In \cite{FT} it is established, actually, that
under the hypotheses of Theorem~\ref{thmFT} one has
$ Q(x, h) \ge c_1 h $ for some $c_1 > 0 $.
As we shall see later, 
using a slight modification of the method of Filaseta and Trifonov, 
we can prove that
$  Q(x, h) \sim \frac{6}{\pi^2} h  $
when 
$   x \to \infty $ and $ \frac{h} {x^{1/5} \log x } \to \infty $.
It is not known whether $ Q(x, h) > 0 $ for smaller $h$, but
in this case we establish 
that $ Q(x, h) $ cannot be much bigger 
than $\frac{6}{\pi^2} h $.

In the present paper we find results of this type
for a more general problem.
Let $ {\bf l} = \langle l_1, \dots, l_r \rangle $ be a vector with 
distinct, non-negative, integer components and
define
\[
  Q_{\bf l}(x) = 
    \sum_{ n \le x} \mu^2(n + l_1) \dots \mu^2(n + l_r) , \qquad
    Q_{\bf l}(x, h) =   Q_{\bf l}(x+h) - Q_{\bf l}(x) .
\]
An asymptotic formula for $ Q_{\bf l}(x) $ in the case
$ {\bf l} = \langle 0, 1 \rangle $ was established in an elementary way
by Carliz \cite{Car}.
Hall \cite{Hall} found an asymptotic formula 
in the general case 
and also proved some results concerning the behavior of $ Q_{\bf l}(x, h) $
on average.
Later Heath-Brown \cite{H-Br} considered again the
particular case $ {\bf l} = \langle 0, 1 \rangle $ and, using his
square sieve, improved the estimate of the error term
in the asymptotic formula for $ Q_{\bf l}(x)$.
Finally, Tsang \cite{Tsang} 
applied the Buchstab--Rosser sieve as well as
Heath-Brown's method and proved that
if 
$ r \le \frac{1}{25} ( \log x / \log \log x ) $
and 
$l_1, \dots, l_r \le c x$ for some constant $c > 0 $,
then
\begin{equation} \label{*l1}
   A({\bf l}) \, x + O \left( r^{12/5} x^{3/5 } (\log x )^{-8/5}  \right) \, \le
     Q_{\bf l}(x)  \, \le
   A({\bf l}) \, x + O \left( r^2 x^{7/11 } (\log x )^7 \right) ,
\end{equation}
where the constants in the $O$-terms depend only on $c$,
\begin{equation} \label{*l2}
  A({\bf l}) = \prod_p \left( 1 - \frac{u(p)}{p^2} \right) 
\end{equation}
(the product is taken over all primes $p$) and where $u(p)$ is the number of
distinct residue classes modulo $p^2$ represented by the integers
$l_1, \dots, l_r$.

Our main result is the following theorem, which states that 
the number of $r$-tuples of squarefree integers,
lying even in a very short interval,
cannot be much bigger than the expected value.
\begin{thm} \label{thm2}
Let $x$, $h$ be real numbers such that $ 10^3 \le h \le x$
and let $r,  l_1, \dots, l_r$ be integers satisfying
\begin{equation} \label{*l7}
  1 \le  r \le \frac{ \log h } { \log \log h } , \qquad
       0 \le l_1 <  \dots < l_r \le x .
\end{equation}
Then we have
\begin{equation} \label{*l8}
  Q_{\bf l} (x, h) \le A({\bf l}) \, h \, \left\{
1 + O \left(  h^{-1/3 \, + \,\rho(h)  }  \right) \right\} ,
    \qquad   \rho(h) = 2 \frac{ \log \log \log h }{ \log \log h } .
\end{equation}
where the constant in the $O$-term is absolute.
\end{thm}

For the proof we apply Selberg's sieve.
We note that the upper bound for $  Q_{\bf l}(x, h) $, 
given by \eqref{*l8}, does not depend on $x$.

Our second result is a generalization of Theorem~\ref{thmFT}.
We apply again sieve methods 
(Buchstab's identity and the sieve of Eratosthenes)
as well as a version of the main proposition of \cite{FT}
and prove the following:
\begin{thm} \label{thm1}
Let $x $ be sufficiently large and $\psi(x)$ be a monotonically increasing function, such that
\begin{equation} \label{*l3}
  2 \le \psi(x) \le e^{-10} (\log x)^{2/3} .
\end{equation}
Suppose that the integers $r, l_1, \dots, l_r$ satisfy
\begin{equation} \label{*l4}
  1 \le r , \qquad
  e^{10 \sqrt{r} } \le 
  (\log x)^{2/3} \psi(x)^{-1} , 
        \qquad   0 \le l_1 <  \dots < l_r \le x 
\end{equation}
and let
\begin{equation} \label{*l5}
  h \ge e^{10 \sqrt{ r } } \, \psi(x) \, x^{1/5} \log x   .
\end{equation}
Then we have
\begin{equation} \label{*l6}
   Q_{\bf l}(x, h) =  A({\bf l}) \, h \, \big\{ 1  + O \left(  \psi(x)^{-1} \right) \big\} ,
\end{equation}
where the constant in the $O$-term is absolute.
\end{thm}

We note that a week version
of Theorem~\ref{thm2} can be deduced from the proof of Theorem~\ref{thm1} 
(see \eqref{*l20}, \eqref{*l25} and \eqref{*l26} ).
More precisely, if the conditions \eqref{*l3} and \eqref{*l4} hold and if 
$ x^{\varepsilon} < h \le x $, where $ \varepsilon > 0 $ is arbitrarily small,
then 
\[
  Q_{\bf l} (x, h) \le A({\bf l}) \, h \, \big\{ 1 + O \left( \psi( h )^{-1} \right)  \big\} 
\]
(the constant in the $O$-term depends on $\varepsilon$).
However the estimate of the remainder term in \eqref{*l8} is much sharper
and also in Theorem~\ref{thm2} we do not impose a lower bound for $h$ depending on $x$.

\paragraph{Acknowledgement:}
The author is grateful to Professor D.R.Heath-Brown
for reading the preliminary version of the paper and for some useful suggestions and remarks.
The author also thanks Plovdiv University Scientific Fund
for the financial support under grant 05-M-45.

\section{Notations, lemmas and some simple estimates.} 

\indent

As usual, $\mu(n)$ is the M\"obius function and $\nu(n)$ denotes the number
of distinct prime factors of $n$. 
The letters $p$ and $q$ are reserved for prime numbers.
By $(k_1, k_2)$ and $[k_1, k_2]$ we denote the greatest common divisor and, respectively, the least common
multiple of the integers $k_1$ and $k_2$.
In this way we also denote open and, respectively, closed intervals, 
but the meaning is always clear from the context.
We write $\# \mathcal M $ for the cardinality of the 
finite set $\mathcal M$.
If it is not specified explicitly, the constants in the $O$-terms and
$\ll$-symbols are absolute.

For a positive integer $k$ we define
\[
 \sigma(k) = \prod_{ p^2 \mid k } p .
\]
If $n$ is a positive integer and if $ {\bf l }= \langle l_1, \dots, l_r \rangle $ is a vector with 
non-negative integer components we define
\[
 \xi(n) = \xi_{\bf l}(n) = \prod_{ j=1 }^r \sigma(n+l_j) 
\]
For any real $z \ge 2$ we denote
\[
    P(z) = \prod_{ p < z } p .
\]

We write $\m D(z, k)$ for an abbreviation of the condition
$ ( \sigma(k), P(z) ) = 1 $.
If $ z_1, z_2, \dots, z_r \ge 2 $ then we introduce another condition
$\m E_{n, {\bf l}} (z_1, z_2, \dots, z_r) $, which means that $n$ satisfies 
$\m D(z_i, n+l_i) $ for all $i=1, \dots, r$.

We shall see that under the conditions
\begin{equation} \label{LL1}
  1 \le h \le x , \qquad 0 \le l_1, \dots, l_r \le x 
\end{equation}
we have
\begin{equation} \label{*l17}
   Q_{\bf l}(x, h) =
      \# \{ n \in (x, x+h] \; : \; \m E_{n, {\bf l}} (\, 2\sqrt{x}, \, 2\sqrt{x}, \dots, 2\sqrt{x} \, ) \; \} .
\end{equation}
Indeed, using \eqref{LL1}, we find that 
$n + l_i < 4x $ 
for any $n \in (x, x+h] $  and for $i = 1, \dots, r$.
Therefore the integer $n+l_i$ is squarefree if and only if 
the condition $\m D(2\sqrt{x}, n+l_i)$ holds.
This implies the representation \eqref{*l17}.

For any squarefree integer $d$ we define
\begin{equation} \label{*l8.5}
   u(d) = \prod_{ p \mid d } u(p)
\end{equation}
and
\begin{equation} \label{*l9}
   N_d(x, h) = \# \{ n \in (x, x+h] \; : \;  \xi(n) \equiv 0 \pmod{ d } \; \} .
\end{equation}

As it is mentioned in \cite{Tsang}, p. 269, 
the congruence $ \xi(n) \equiv 0 \pmod{p} $ has $u(p)$ solutions modulo $p^2$
and, respectively, the congruence $ \xi(n) \equiv 0 \pmod{ d } $,
where $d$ is squarefree,  has exactly
$ u(d) $ solutions modulo $ d^2 $.
Therefore we find
\begin{equation} \label{*l10}
    N_{d}(x, h) = h \frac{u(d)}{d^2} + O(u(d)) .
\end{equation}

Obviously, for any prime $p$ we have
\begin{equation} \label{*l11}
u(p) \le r .
\end{equation}
We can also assume that
\begin{equation} \label{*l12}
u(p) \le p^2 - 1 
\end{equation}
for all $p$ because otherwise we would have $ Q_{\bf l} (x, h) = A({\bf l}) = 0 $ 
and our results would be trivial.

We shall use the following simple estimate:
\begin{equation} \label{*l13}
 1 \le  A({\bf l})^{-1} \le    e^{ 9 \sqrt{r} }   .
\end{equation} 
Indeed, the first inequality is obvious. To prove the second one we 
use \eqref{*l2} to write
\begin{equation} \label{*l14}
 A({\bf l})^{-1} = 
    \prod_{  p \le \sqrt{2r} } \left( 1 - \frac{u(p)}{p^2} \right)^{-1} \; 
    \prod_{  p > \sqrt{2r} } \left( 1 - \frac{u(p)}{p^2} \right)^{-1} =
P_1 P_2 ,
\end{equation}
say. 
From \eqref{*l12} and from 
the well-known upper bound 
in the Tchebishev prime number theorem
\[
   \prod_{p \le w} p \le 4^w \qquad \text{ for } \qquad w \ge 1
\]
we find
\begin{equation} \label{*l15}
 P_1 \le 
\prod_{  p \le \sqrt{2r} } \left( 1 - \frac{p^2-1}{p^2} \right)^{-1} =
\prod_{  p \le \sqrt{2r} } p^2 \le 4^{ 2 \sqrt{2 r } } .
\end{equation}
Respectively, from \eqref{*l11}
we get
\begin{equation} \label{*l16}
  \log P_2 \le 
      - \sum_{  p > \sqrt{2r} } \log ( 1 - r p^{-2} )
   \le 2 r \sum_{ n > \sqrt{2r} } n^{-2}  \le 
        2 \sqrt{2r}  .
\end{equation}
The right inequality in \eqref{*l13} follows from
\eqref{*l14} -- \eqref{*l16}.

\bigskip

The core of the proof of Filaseta and Trifonov's theorem is 
the following:
\begin{lem}[Filaseta, Trifonov] \label{lemFT}
Suppose that $x$ is sufficiently large and $h = c x^{1/5} \log x$, 
where $c>1$ is sufficiently large constant.
Then there exists a constant $ \gamma > 0 $ such that
\[
  \#  \{ d \in ( \, h\sqrt{\log x}, 2 \sqrt{x} \, ] \cap \mathbb Z \; : \;
            m d^{ \, 2} \in (x, x+h] \;\, \text{ for some } \;\, m \in \mathbb Z \; \} 
    \ll c^{-\gamma} h .
\]
\end{lem}
The proof of this result is presented in detail in \cite{FT}.
To establish our Theorem~\ref{thm1} we shall use the following modification:
\begin{lem} \label{lem1}
Suppose that $x$ is sufficiently large and
\begin{align}
 &
 x \le X \le 2x , \qquad \qquad
 1 \le R \le (\log x)^{2/3} , 
      \notag \\
 &
  h \ge R \, x^{1/5} \, \log x ,
 \qquad \;\,
 \lambda = R^{ -1} \, h \, \log x \le 2\sqrt{x} .
   \notag
\end{align}
Then we have
\[ 
   \# \{ d \in [ \, \lambda, 2\sqrt{x} \, ] \cap \mathbb Z \; : \;
            m d^{ \, 2} \in (X, X+h] \;\, \text{ for some } \;\, m \in \mathbb Z \; \} 
                  \;  \ll  \; R^{-1}  h .
\]
\end{lem}
The proof differs very slightly from the proof
of Lemma \ref{lemFT}, so we omit it.

\section{ Proof of Theorem \ref{thm2}.}

\indent

We may assume that $ h \ge h_0$, where $h_0$ is a sufficiently large absolute constant.
We use the representation of 
$  Q_{\bf l}(x, h) $ in the form \eqref{*l17} and,
since the condition $ \mathcal E_{ n, {\bf l} }( \, 2\sqrt{x}, \, 2\sqrt{x}, \dots, 2\sqrt{x} \, ) $ 
is equivalent to $ (\xi(n), P(2\sqrt{x}) ) = 1 $, we can write
\[
  Q_{\bf l}(x, h) = \sum_{ x < n \le x+h } \; \sum_{ d \mid (\xi(n), P(2\sqrt{x}) ) } \mu(d) .
\]

Now we apply Selberg's upper bound sieve.
Let $\lambda(d)$ be real numbers defined for squarefree integers $d$.
We suppose that 
$\lambda(1) = 1 $ and $ \lambda(d) = 0 $ for $d > z$,
where $z$ is a parameter for which we assume
\begin{equation} \label{*l31}
 2 < z < 2\sqrt{x} .
\end{equation}
We note that if $d$ is squarefree and $d \le z$ then $d \mid P(2\sqrt{x}) $.
Hence we find
\begin{align}
  Q_{ \bf l }(x, h) 
     & \le 
           \sum_{ x < n \le x+h } \left( \sum_{ d \mid \xi(n)  } \lambda(d) \right)^2 
            = \sum_{ x < n \le x+h }  \; \sum_{ d_1, d_2 \mid \xi(n) } \lambda(d_1) \lambda(d_2) 
                  \label{*l32}  \\
     & \notag \\
     & = \sum_{ d_1, d_2 \le z } \lambda(d_1) \lambda(d_2) \;
                   N_{[d_1, d_2]}(x, h) ,
            \notag 
\end{align}
where $  N_d(x, h) $ is defined by \eqref{*l9}.

We use \eqref{*l10} and \eqref{*l32} to get
\begin{equation} \label{*l33}
  Q_{\bf l}(x, h) \le h V + O(R \,) ,
\end{equation}
where
\begin{equation} \label{*l34}
   V = \sum_{ d_1, d_2 \le z } \frac{ \lambda(d_1) \lambda(d_2) } { [d_1, d_2]^2 } \, u([d_1, d_2]),
   \; \,
   R = \sum_{ d_1, d_2 \le z } |\lambda(d_1)| \, |\lambda(d_2)| \, u([d_1, d_2]).
\end{equation}

We define $\lambda(d)$ for $1 < d \le z $ 
in such a way as to minimize $V$. 
By a straightforward application of Selberg's method
(see, for example, \cite{HR}, Chapter 3) we can verify that the optimal choice is
\begin{equation} \label{*l35}
  \lambda(d) =  \mu(d)  \prod_{ p \mid d } \, \left( 1 - \frac{u(p)}{p^2} \right)^{-1} \,
                \frac{ H(z/d, d) }{ H(z) } ,
\end{equation}
where
\begin{equation} \label{*l36}
  H(y, m) = \sum_{ \substack{ k \le y \\ (k, m)=1 }} \frac{ \mu^2(k) \, u(k) } { k^2 } \prod_{ p \mid k }
                          \left( 1 - \frac{u(p)}{p^2} \right)^{-1} , 
  \qquad H(y) = H(y, 1) .
\end{equation}
In this case the minimal value of $V$ is
\begin{equation} \label{*l37}
   V_{ \min } = H(z)^{-1} .
\end{equation}
We leave the calculations to the reader.

From \eqref{*l2}, \eqref{*l35} and \eqref{*l36} we easily find
\begin{equation} \label{*l38}
   |\lambda(d)| \le \mu^2(d) \, \prod_{ p \mid d } \left( 1 - \frac{u(p)}{p^2} \right)^{-1}
     \le  A({\bf l})^{-1} \mu^2(d).
\end{equation}

Having in mind \eqref{*l33}, \eqref{*l34} and \eqref{*l37} we obtain
\begin{equation} \label{*l39}
  Q_{\bf l}(x, h) \le h H(z)^{-1} + O(G^{ \, 2} ) ,
\end{equation}
where
\[
    G = \sum_{ d \le z } |\lambda(d)| u(d) .
\]
From \eqref{*l8.5}, \eqref{*l11} and \eqref{*l38}  we find
\begin{equation} \label{*l40}
  G  \le  A({\bf l})^{-1} \sum_{ d \le z } r^{\nu(d)} \mu^2(d)
        = A({\bf l})^{-1} U ,
\end{equation}
say.
To estimate the sum $U$ we assume that 
\begin{equation} \label{velikden}
 \nu = 1 + \frac{r}{\log z} \le 2
\end{equation}
and, using Euler's identity and the elementary properties of Riemann's zeta-function, we easily get
\begin{align} 
  U 
     & \le 
       \sum_{ d=1 }^{\infty} \left( \frac{z}{d} \right)^{\nu} \mu^2(d) \, r^{\nu(d)}
        = z^{\nu} \prod_{ p } \left( 1 + \frac{r}{p^{\nu} } \right) 
           \le z^{\nu} \left( \frac{\nu}{\nu-1} \right)^r
              \label{*l41} \\ 
     & \le z \, ( 2 \, e \, r^{-1} \log z )^{r} .
               \notag
\end{align}

Applying Euler's identity we also find
\[
   \sum_{ k = 1}^{\infty} \frac{ \mu^2(k) \, u(k) } { k^2 } \prod_{ p \mid k }
                          \left( 1 - \frac{u(p)}{p^2} \right)^{-1}  = A({\bf l})^{-1} .
\]
Hence, if we define $\omega(z)$ by 
\begin{equation} \label{*l42}
  H(z) =  A({\bf l})^{-1} - \omega(z) ,
\end{equation}
then, using \eqref{*l2}, \eqref{*l8.5} and \eqref{*l11} we get
\begin{align}
    0 \le  \omega(z) 
     & =  
         \sum_{ k > z } \frac{ \mu^2(k) u(k) }{ k^2} 
                       \prod_{ p \mid k } \left( 1 - \frac{u(p)}{p^2} \right)^{-1}
           \label{*l43} \\
     &
      \le
         A({\bf l})^{-1} \sum_{ k > z } \frac{ \mu^2(k) \, r^{\nu(k)} }{ k^2} 
     =  A({\bf l})^{-1} U_1 ,
       \notag 
\end{align}
say.
Arguing as in the proof of \eqref{*l41} we find
\[
 U_1 \ll z^{-1} ( 2 \, e \, r^{-1} \log z )^{r} .
\]
The last estimate, \eqref{*l42} and \eqref{*l43} imply
\begin{equation} \label{*l44}
  H(z)^{-1} =  A({\bf l}) \, \big\{ 1 + O ( z^{-1} ( 2 \, e \, r^{-1} \log z )^{r}  ) \big\} .
\end{equation}

From \eqref{*l39}, \eqref{*l40}, \eqref{*l41} and \eqref{*l44} we obtain
\begin{equation} \label{*l45} 
  Q_{\bf l } (x, h) \le  A({\bf l}) \, h \, \big\{ 1 + O ( \Delta )  \big\} ,
\end{equation}
where
\[
   \Delta =  
     z^{-1} ( 2 e \, r^{-1} \log z)^{r}   +  A({\bf l})^{-3} h^{-1} z^2 ( 2 e \, r^{-1} \log z)^{2r} .
\]
We choose 
\[
z = h^{1/3} \, ( r^{-1} \log h )^{-r/3} .
\]
and note that
the conditions \eqref{*l31} and \eqref{velikden} are satisfied.
Using  \eqref{*l7} and \eqref{*l13} we find
\[
  \Delta \ll h^{-1/3 + \rho(h) } ,
\]
where $\rho(h)$ is specified by \eqref{*l8}.
We leave the standard calculations to the reader,

From the last formula and \eqref{*l45} we obtain \eqref{*l8},
so Theorem~\ref{thm2} is proved.

\section{ Proof of Theorem \ref{thm1}.}

\indent

We can assume that $h \le x^{7/11} (\log x)^{10} $ 
because otherwise \eqref{*l6} is a consequence of \eqref{*l1}, \eqref{*l4} and \eqref{*l13}.
Suppose that $ \mathcal A \subset (x, x+h] \cap \mathbb Z $ and let
\begin{equation} \label{*l18}
   \lambda_0 = e^{10 \sqrt{r} } \, \psi(x) .
\end{equation}
Using the Buchstab identity 
(see, for example, \cite{HR}, Chapter 7, p. 204),
we find for any $i=1, \dots, r$ that
\begin{align}
   &
      \# \{ n \in \mathcal A \; : \; \m D(2\sqrt{x}, n+l_i) \; \} 
       \label{*l19} \\
   & \qquad \qquad \qquad \qquad  
       =  \# \{ n \in \mathcal A \; : \; \m D(\lambda_0, n+l_i) \; \}  
            \notag \\
     & \qquad \qquad \qquad  \qquad \;
      -  \sum_{ \lambda_0 \le q < 2\sqrt{x} } 
         \# \{ n \in \mathcal A \; : \; q^2 \mid n+l_i \; , \; \m D(q, n+l_i)  \; \} .
      \notag
\end{align}

Having written $ Q_{\bf l}(x, h) $ 
in the form \eqref{*l17},
we apply \eqref{*l19} with $i=1$ and with the set $\m A$ consisting of the integers
$n \in (x, x+h] $ satisfying the conditions
$\m D(2\sqrt{x}, n+l_j) $ for all $j=2, \dots, r$.
We get
\begin{align}
 Q_{\bf l}(x, h) 
      & =
      \#  \{ n \in (x, x+h] \; : \; \m E_{n, {\bf l}} (\lambda_0, 2\sqrt{x}, \dots, 2\sqrt{x} \, ) \; \} 
               \notag \\
      & \;
       -  \sum_{ \lambda_0 \le q < 2\sqrt{x} } 
        \# \{ n \in (x, x+h]  \; : \; q^2 \mid n+l_1 \; , \; 
            \m E_{n, {\bf l}} (q, 2\sqrt{x}, \dots, 2\sqrt{x} \, ) \; \} .
                \notag 
\end{align}

We consider the first term from the right side of the last identity
and apply \eqref{*l19} again, this time with $i=2$ and with the set 
$ \mathcal A $ consisting of all integers $n \in (x, x+h]$, which satisfy
$\m D(\lambda_0, n+l_1) $ and
$\m D(2\sqrt{x}, n+l_j) $ for all $j=3, \dots, r $.
In this way we find another identity for  $ Q_{\bf l}(x, h) $.

Proceeding in this manner we obtain
\begin{equation} \label{*l20}
 Q_{\bf l}(x, h) =
   R_0 - \sum_{ \nu = 1 }^r \sum_{ \lambda_0 \le q < 2\sqrt{x} } R_{ \nu, q} =
   R_0 - \Sigma , 
\end{equation} 
say, where
\begin{align} 
    R_0 
    & = 
      \# \{ n \in (x, x+h] \; : \; \m E_{ n, \bf l} (\lambda_0, \lambda_0, \dots, \lambda_0) \; \} ,
      \notag \\
   R_{\nu, q} 
     & = 
     \# \big\{ n \in (x, x+h] \; : \; q^2 \mid n+l_{\nu} \; , \;
           \m E_{ n, \bf l} \big(
\underbrace{\lambda_0, \dots, \lambda_0}_{ \nu -1 }, 
\, q, \, 
\underbrace{ 2\sqrt{x}, \dots, 2\sqrt{x} }_{ r - \nu }  \, \big) \; \big\} .
              \notag 
\end{align}

Consider $R_0$. 
It is clear that the condition 
$ \m E_{ n, \bf l} (\lambda_0, \lambda_0, \dots, \lambda_0) $
is equivalent to
$ (\xi(n), P(\lambda_0) ) = 1 $.
Therefore
\begin{equation} \label{*l21}
  R_0 = 
     \sum_{ x < n \le x+h } \; \sum_{ d \mid ( \xi(n), P(\lambda_0) ) } \mu(d) =
       \sum_{ d \mid P(\lambda_0) } \mu(d) 
                N_d(x, h) ,
\end{equation}
where $ N_d(x, h) $ is defined by \eqref{*l9}.
Using \eqref{*l10} and \eqref{*l21} we get
\begin{equation} \label{*l22}
    R_0 = h W + O(H) ,
\end{equation}
where
\[
 W = \sum_{ d \mid P(\lambda_0) } \mu(d) \frac{u(d)}{d^2}  = 
   \prod_{ p < \lambda_0 } \left( 1 - \frac{u(p)}{p^2} \right) , \qquad
   H = \sum_{ d \mid P(\lambda_0) } u(d) .
\]

Consider $W$ and $H$.
Arguing as in the proof of \eqref{*l13}
and using \eqref{*l2}, \eqref{*l11} and \eqref{*l18} we find
\begin{equation} \label{*l23}
   W =
      \prod_{ p  } \left( 1 - \frac{u(p)}{p^2} \right) \;
    \prod_{ p \ge \lambda_0 } \left( 1 - \frac{u(p)}{p^2} \right)^{-1}
      =  A({\bf l}) \, \left( 1 + O \left( \frac{r}{\lambda_0} \right) \right) .
\end{equation}
Respectively, we have
\begin{equation} \label{*l24}
   H \le \sum_{ d \mid P(\lambda_0) } r^{\nu(d)} = \prod_{ p < \lambda_0 } (1+r) \le (1+r)^{\lambda_0} .
\end{equation}

From \eqref{*l22} -- \eqref{*l24} we get
\[ 
   R_0 =  A({\bf l}) \, h \, \left( 1 + O ( \Delta )  \right)  ,
\]
where
\[
  \Delta = r \lambda_0^{-1} + h^{-1}  A({\bf l})^{-1} (1+r)^{\lambda_0} .
\]
Ii is not difficult to verify,
using \eqref{*l4}, \eqref{*l5}, \eqref{*l13} and \eqref{*l18},
that
$ \Delta \ll \psi(x)^{-1} $
and we find
\begin{equation} \label{*l25}
     R_0 =  A({\bf l}) \, h \, \left( 1 + O ( \psi(x)^{-1} \right)  .
\end{equation}

Consider now the sum $ \Sigma $, specified by \eqref{*l20}.
We have
\begin{equation} \label{*l26}
   0 \le   \Sigma  \le r \, \max_{ 1 \le \nu \le r } S_{ \nu} .
\end{equation}
where
\[
   S_{\nu} = \sum_{ \lambda_0 \le q < 2\sqrt{x} } 
          \# \{ n \in (x, x+h] \; : \; q^2 \mid n+l_{\nu} \; \} .
\]
We shall prove that
\begin{equation} \label{*l30.5}
  S_{\nu}  \ll  e^{-10 \sqrt{r} } \, \psi(x)^{-1} \, h .
\end{equation}

Define
\begin{equation} \label{*l28}
 \lambda =   e^{-10 \sqrt{r} } \, \psi(x)^{-1}  \, h \, \log x  .
\end{equation}
First we assume that $ \lambda < 2\sqrt{x} $. 
In this case we divide $ S_{\nu} $ into two parts
\begin{equation} \label{*l27}
  S_{\nu}  = S' + S'' .
\end{equation}
In $S'$ the summation is taken over the primes $q \in [ \lambda_0,  \lambda  ) $
and $S''$ is the contribution from the primes 
$q \in [ \lambda, 2\sqrt{x}) $.

Using Tchebishev's prime number theorem and  
\eqref{*l4}, \eqref{*l5}, \eqref{*l18} and \eqref{*l28}
we get
\begin{align} 
    S'
     & =
      \sum_{ \lambda_0 \le q < \lambda } 
        \sum_{ \substack{ x < n \le x+h \\ n + l_{\nu} \equiv 0 \pmod{q^2} }} 1
           \le
         \sum_{ \lambda_0 \le q < \lambda }  \left( \frac{h}{q^2} + 1 \right)
\ll
         \frac{h}{\lambda_0} +  \frac{\lambda}{\log \lambda} 
        \label{*l29} \\
       & \notag \\
         &   \ll e^{-10 \sqrt{r} } \, \psi(x)^{-1} \, h .
             \notag 
\end{align}
To estimate $S''$ we write it in the form
\[
   S'' = \sum_{ \lambda \le q < 2\sqrt{x} } 
          \# \{ m \in (x+l_{\nu}, x+l_{\nu}+h] \; : \; q^2 \mid m \; \} 
\]
and apply Lemma \ref{lem1} with $X = x+l_{\nu}$ and $ R = e^{10 \sqrt{r} } \psi(x) $.
We find
\begin{equation} \label{*l30}
  S''  \ll  e^{-10 \sqrt{r} } \, \psi(x)^{-1} \, h .
\end{equation}
The estimate \eqref{*l30.5} follows from \eqref{*l27} -- \eqref{*l30}.

In the case $ \lambda \ge 2\sqrt{x} $ we proceed as in the estimation of 
$S'$ and easily find that \eqref{*l30.5} holds as well.

From \eqref{*l13}, \eqref{*l20}, \eqref{*l25} -- \eqref{*l30.5}
and \eqref{*l27}  we obtain \eqref{*l6}
and Theorem~\ref{thm1} is proved.

\bigskip
\bigskip

\vbox{
\hbox{Department of Mathematics}
\hbox{Plovdiv University ``P. Hilendarski"}
\hbox{24 ``Tsar Asen" str.}          
\hbox{Plovdiv 4000}
\hbox{Bulgaria}
\smallskip
\hbox{Email: dtolev@pu.acad.bg}}

\end{document}